\input amstex
\documentstyle{amsppt}\magnification=\magstep1
\magnification=1200
%\hcorrection{0.25in}
%\vcorrection{-0.8in}
\NoRunningHeads \NoBlackBoxes \baselineskip=12pt
\vsize=8.5truein \hsize=6.0truein \hoffset=0.3truein
\voffset=0.2truein \loadmsam \loadbold \loadmsbm \loadeufm
%\loadeusm \UseAMSsymbols
%\font\boldtitlefont=cmb10 scaled\magstep2 \font\sectionfont=cmb10
%scaled\magstep1
 \NoRunningHeads
 %\nologo
\document
\define\({\left(}
\define\){\right)} \define\[{\left[} \define\]{\right]}
 
\define\aint{\rlap{\kern2pt\vrule height2.5pt width9.5pt depth-2.1pt}
\int_{B_{\rho}}} \TagsOnRight
  % \nologo

\topmatter
\title
The smoothness of Riemannian submersions with nonnegative
sectional curvature
\endtitle
\author   Jianguo Cao and     Mei-Chi Shaw  \endauthor
 \address Jianguo Cao,
 Department of Mathematics, University of Notre Dame,
Notre Dame, IN  46556 USA. \endaddress \email
cao.7$\@$nd.edu\endemail
\address Mei-Chi Shaw,
Department of Mathematics, University of Notre Dame, Notre Dame,
IN  46556 USA  \endaddress \email mei-chi.shaw.1\@nd.edu \endemail

\thanks Both authors are supported   by NSF Grants. Research at MSRI is
supported in part by NSF grant DMS-9810361. The first author is
grateful to MSRI and Max-Planck Institute for Mathematics at
Leipzig for their hospitality.
\endthanks
\endtopmatter
\bigskip

   In this article, we study the smoothness of Riemannian
   submersions for open manifolds with non-negative sectional
   curvature.
    Suppose that $M^n$ is a $C^\infty$-smooth,
   complete and non-compact Riemannian manifold with nonnegative
   sectional curvature. Cheeger-Gromoll [ChG] established a
   fundamental theory for such a manifold. Among other things,
   they showed that $M^n$ admits a totally convex exhaustion $\{
   \Omega_u\}_{u \ge 0}$ of $M^n$,  where $\Omega_0 = \Cal S$ is a
   totally geodesic and  compact submanifold without boundary. Furthermore,
   $M^n$ is diffeomorphic to the normal vector bundle of the soul
   $\Cal S$.

       Sharafutdinov found  that there exists a distance
       non-increasing retraction $\Psi: M^n \to \Cal S$ from the
       open manifold $M^n$ of non-negative sectional curvature to
       its soul, (cf. [Sh], [Y2]).  Perelman [Per] further showed that such a map
       $\Psi$ is indeed a $C^1$-smooth Riemannian submersion.
       Furthermore, $\Psi[Exp_q(t \vec v)] = q$ for any $q \in
       \Cal S$ and $\vec v \bot T_q(\Cal S)$. Therefore, the fiber
       $F_q = \Psi^{-1}(q)$ is a $k$-dimensional submanifold,
       which is $C^\infty$-smooth almost
       everywhere, where $k = dim(M^n) -dim(\Cal S)>0$.

       Guijarro [Gu]  proved that the fiber $F_q$ is indeed a
       $C^2$-smooth submanifold for each $q\in S$. In this paper,
       we prove that the fibres are $C^\infty$-smooth.

\proclaim{Theorem 1}  Let $M^n$ be a complete, non-compact and
$C^\infty$-smooth  Riemannian manifold with nonnegative sectional
curvature. Suppose $\Cal  S$ is a soul of $M^n$. Then any distance
non-increasing retraction $\Psi: M^n \to \Cal S$ must give rise to
a $C^\infty$-smooth Riemannian submersion.

 Consequently, if $\Bbb R^k = \Cal N_q(\Cal S, M^n)$ is the normal
 space of the soul $\Cal S$  in $M^n$ at $q$, then the fiber $F_q = \Psi^{-1}(q) =
 \text{Exp}_q (\Bbb R^k)$ is a $k$-dimensional $C^\infty$-smooth
 submanifold of $M^n$, for any $q \in S$.
\endproclaim

Professor Wilking kindly informed us that he
      has  recently obtained a similar result
        (cf. [Wi]).
        His method is completely
       independent of ours.
   Our  proof of Theorem 1
   uses a flat strip theorem associated with Cheeger-Gromoll exhaustion
   (cf. Theorem 4 below),  an uniform estimate for cut-radii of  convex
   subsets in [ChG] and a smooth extension theorem for  ruled surfaces.

    For each compact convex subset $\Omega \subset M^n$, we
    let $U_\epsilon (\Omega) = \{ x \in M^n | d(x, \Omega) < \epsilon \}$.
    Its cut-radius is given  by $\delta_\Omega = \sup\{\epsilon |
    \text{ there is a unique nearest point projection  } \Cal P_\Omega:
    U_\epsilon (\Omega) \to
    \Omega   \}$.

        For each $x \in M^n$, we let $\text{Inj}_{M^n}(x)$ be the
        injectivity radius of $M^n$ at $x$. Similarly, let
        $\text{Inj}_{M^n}(A) = \sup\{ \text{Inj}_{M^n}(x)| x \in A
        \}$.

A subset $\Omega$ of a complete Riemannian manifold $M^n$ is said
to be {\it totally convex} if for any pair of points $\{p, q\}
\subset \Omega$ and for any geodesic segment $\sigma$ joining $p$
and $q$, the geodesic segment $\sigma$ is contained in $\Omega$.
There is a totally convex exhaustion $\{
   \Omega_u\}_{u \ge 0}$ of $M^n$ given in [ChG]. By comparing the inner
angles of geodesic triangles, we have the following semi-global estimate for cut-radius.

     \medskip
     \proclaim{Lemma 2}(Lemma 2.4 of [ChG], [CaS])  Let $A \subset
\Omega_{T}$ be a connected, convex and compact subset  in a
  Riemannian manifold $M^n$ with nonnegative curvature, let
  $K_0 =\max\{K(x) | x \in \Omega_{T+1} \}$, $\text{Inj}_{M^n}(\Omega_{T})$ be the upper bound of sectional curvature
  on $ \Omega_{T+1} $ and $\Cal S$ be as above. Suppose that
$\dim(\Omega_{T}) =n$. Then the subset $A$ has cut-radius bounded
below by
$$
\delta_A \ge \delta_0(T) = \frac{1}{4}\min\{
\text{Inj}_{M^n}(\Omega_{T}), \frac{\pi}{\sqrt{K_0}}, 1 \},
$$
where $\delta_{0}(T)$ is independent of choices of  $A$ with $A
\subset \Omega_{T}$.
\endproclaim

Let us briefly recall the Cheeger-Gromoll convex exhaustion.
 According to [ChG], there is a partition $a_0 =0 < a_1 < ... a_m < a_{m+1} = \infty$ of
 $[0, \infty)$ and an exhaustion $\{\Omega_u\}_{u \ge 0}$ of
 $M^n$ such that the following holds:
 \medskip \noindent
 (1) $M^n = \cup_{u \ge 0}\Omega_u$. If $u > a_m$ then  $dim[\Omega_u] = n $.
  If  $u \le a_m$, then $\dim[ \Omega_{u}  ]<
 n$.
\medskip \noindent
 (2) $\Omega_0 = \Cal S$ is the soul of $M^n$, which  is a totally
 geodesic $C^\infty$-smooth compact submanifold without boundary.
 \medskip \noindent
 (3) If $u > 0$,  $\Omega_u$ is a
 totally convex, compact subset of $M^n$ and hence $\Omega_u$  is
 a compact
 submanifold with a $C^\infty$-smooth relative interior. Furthermore,
 $\dim(\Omega_u) = k_u > 0$ and
 $\Omega_u$ has a non-empty $(k_u -1)$-dimensional relative boundary
 $\partial \Omega_u$;
\medskip \noindent
(4) For any $u_0 \in [a_{j}, a_{j+1}]$ and $0 \le t \le u_0 -
a_{j}$, the family $\{ \Omega_{u_0 -t}\}_{t \in [0, u_0 - a_{j}]}$
is given by {\it the inward equidistant evolution}:
$$
  \Omega_{u_0-t} = \{ x \in \Omega_{u_0 }| d(x, \partial \Omega_{u_0}) \ge t
  \}.  \tag2.1
$$
\medskip \noindent
(5) If $u > a_m $ then $u-a_m = \max\{d(x, \partial \Omega_u) | x
\in \Omega_{u}) \}$. If $0 \le j \le m-1$ then $a_{j+1} - a_{j} =
\max\{d(x, \partial \Omega_{a_{j+1}}) | x \in  \Omega_{a_{j+1}}
\}$ and hence $\dim[ \Omega_{a_{j}} ] < \dim [\Omega_{a_{j+1}} ]$
for $j \ge 0$.

\medskip

   Assume that $k = \dim[M^n] - \dim[\Cal S ] = \dim (F_q) $ for
   all $q \in \Cal S$.
   Since $M^n = \cup_{T \ge 0}\Omega_T$, it is sufficient to verify
   that the subset
   $[U_{\delta_0(T)}(\Omega_T) \cap F_q]$ has a $k$-dimensional
   $C^\infty$-smooth interior, where $\delta_0(T)$ is given by
   Lemma 2 and $T > a_m$.

    For this purpose, we need to study the geometry of
    the equidistant hypersurfaces from $\partial \Omega_u$.
    Federer [Fe] has studied the smoothness of {\it the outward  equidistant
    hypersurfaces} $ \partial [U_\epsilon
    (\Omega)] $ for $0 < \epsilon < \delta_\Omega$. Following his
    approach, we consider the outward normal cone of $\Omega$ as
    follows:
    $$
          \Cal N^+(\Omega, M^m) = \{ (p, \vec v) | p \in \Omega,
       d(Exp_p(t\vec v), \Omega) = t |\vec v|, \text{ for } 0 \le t |\vec v| < \delta_\Omega
          \}.
    $$

If $\{ \Omega_u\}$ is the Cheeger-Gromoll convex exhaustion as
above and $u > 0$, then the relative boundary $\partial \Omega_u$
is not necessarily smooth. We are going to study the corresponding
decomposition of $ \Cal N^+(\Omega, M^m)$:
$$\Cal N^+_{p}(\Omega_u, M^n) \subset \big[\Cal
N^+_{p}(\Omega_u, \text{int}(\Omega_{u+ \epsilon}) ) \oplus \Cal
N^+_{p}(\text{int}(\Omega_{u+ \epsilon}), M^n)\big], \tag2.2
$$
where  $\Cal N^+_{p}(\Omega_u, \text{int}(\Omega_{u+ \epsilon}))$
is  defined by
$$
\aligned
 \Cal N^+\big(\Omega_u, \text{int}(\Omega_{u+ \epsilon}) \big)=  & \{ (p, \vec v) |
  p \in \Omega_u,
       d(Exp_p(t\vec v), \Omega_u) = t |\vec v|, \\
       & \text{ for } 0 \le t |\vec v| < \delta_{\Omega_u},  Exp_p(t\vec v) \in
       \text{int}(\Omega_{u+ \epsilon}) \}.
       \endaligned
       $$

Our next step is to choose $\epsilon $ sufficiently small so that
(1)  there is a nearest point projection  $ \Cal P:
\text{int}(\Omega_{u+ \epsilon}) \to \Omega_u$;  and (2) $\Omega_u
= \{ x \in \Omega_{u+ \epsilon}| d(x, \partial\Omega_{u +
\epsilon}) \ge \epsilon \}$ holds. We first find $j$ so that $
a_j \le u < a_{j +1}$ for some $ 0 \le j \le m$. Let $T = u + a_m
+ 1$ and $\delta_0(T)$ be given by Lemma 2. It follows from a
result of Yim  that there is
        a constant $C_T$ such that,  for $0 \le u_1 < u_2 \le T$,  we have
        $$
          \max\{ d(x, \Omega_{u_1}) | x \in \Omega_{u_2} \} \le
          C_T (u_2 -u_1), \tag2.3
        $$
see [Y2, Theorem A.5(3)]. In what follows, we always choose
    $$ 0 < \epsilon = \epsilon_u <  \min\{[a_{j+1}
- u], \frac{\delta_0(T)}{2C_T}\}, \tag2.4
  $$ where $u \in [a_j, a_{j+1})$, $T= u + a_m +1$ and $\delta_0(T)$ is given by Lemma
  2.

With such a choice of $ \epsilon = \epsilon_u$ by (2.4), the
geometry of $\Cal N^+_p(\Omega_u, \text{int}(\Omega_{u+
\epsilon}))$ is determined by its minimal normal vectors which we
now describe.

\definition{Definition 3 } (Minimal normal vector) Let $\Omega_u$,
$ \Omega_{u+\epsilon}$ and $\Cal N^+(\Omega_u,
\text{int}(\Omega_{u+ \epsilon}))$ be as above. Let  $\sigma_{(p,
\vec v)}: [ 0, \epsilon] \to M^n$ be a geodesic given by $
\sigma_{(p, \vec v)}(t) = \text{Exp}_p(t \frac{\vec v}{|\vec v |})
$, where $\vec v \neq 0$. If $\sigma_{(p, \vec v)} $ is a
length-minimizing geodesic from $p \in \Omega_u$ to $\partial
\Omega_{u+\epsilon}$, then $\vec v$ is called a minimal normal
vector in $\Cal N^+_p(\Omega_u, \text{int}(\Omega_{u+
\epsilon}))$.
\enddefinition

\smallskip

   It is known that any other normal vector $\vec w \in \Cal N^+_p(\Omega_u,
\text{int}(\Omega_{u+ \epsilon}))$ can be expressed as a linear
combination of {\it minimal normal vectors} at $p$. Moreover, the
convex hull of minimal normal vectors at $p$ is equal to $\Cal
N^+_p(\Omega_u, \text{int}(\Omega_{u+ \epsilon}))$, (cf. [Y1,
Proposition 1.7]).

\medskip

For each $p \in M^n$, we let $\Cal V_p = T_p(F_{\Psi(p)})$ and
$\Cal H_p = [\Cal V_p]^\bot$. A geodesic $\alpha: [a, b] \to M^n$
is said to be horizontal, if $\alpha'(t) \bot F_{\Psi(\alpha(t))}$
for all $t \in [a, b]$. We need the following flat strip theorem
for the proof of Theorem 1.

\bigskip
     \proclaim{Theorem 4}   Let $\{\Omega_u\}$ be the Cheeger-Gromoll totally convex
     exhaustion of $M^n$ as above.  Suppose
that $\Psi: M^n \to \Cal S$ be a distance non-increasing
retraction and $F_q = \Psi^{-1}(q)$ be a fibre for some $q \in
\Cal S$. Then for  $p \in F_q \cap \Omega_{u}$ and any $(p, \vec
v) \in \Cal N^+(\Omega_{u}, M^n)$, we have
$$
   \Psi[ Exp_{p}([ \Bbb R\{\vec v\}])] = \Psi(p) = q.  \tag4.1
$$
Moreover, if $\dim(\Cal S) \ge 1$ and if $\vec w \in \Cal H_p$ has
$|\vec w| = 1 = |\vec v|$, then the surface $\Sigma^2_{\vec v,
\vec w} = Exp_p[ \Bbb R \{\vec v \} \oplus  \Bbb R \{\vec w \}]$
is totally geodesic immersed flat plane in $M^n$.
\endproclaim

A result similar to Theorem 4 was proved in [CaS] via a totally
different method.

 \demo{Proof of Theorem 4} Theorem 4 was proved by Perelman [Per] for the
case of $\Omega_0 = \Cal S$. Applying Perelman's argument for the
case of  $p \notin \Cal S$, Guijarro [Gu1] found the following
sufficient condition for (4.1).

\medskip
\noindent (4.2) {\it  $\vec v \in \Cal V_p$ stays vertical under
parallel transport along any horizontal broken geodesic.}

\smallskip

  Guijarro showed that (4.1) follows from (4.2). Moreover, if (4.2)
holds and if $\vec w \in \Cal H_p$ has $|\vec w| = 1 = |\vec v|$,
then the surface $\Sigma^2_{\vec v, \vec w} = Exp_p[ \Bbb R \{\vec
v \} \oplus \Bbb R \{\vec w \}]$ is totally geodesic immersed flat
plane in $M^n$, (cf. Theorem 3.1 of [Gu1]).

\medskip

    In order to see that $\Cal N^+(\Omega_u, M^n) \subset \Cal V_p$ holds, we
     recall that any horizontal geodesic $\alpha$ is contained a tubular
     neighborhood of the soul $\Cal S$, by Perelman's theorem
     [Per]. Hence, $\alpha$ is contained in a compact totally
     geodesic subset $\Omega_T$ for a sufficiently large $T$. It
     follows from Theorem 5.1 of [ChG] that $\alpha \subset
     \partial \Omega_u$ for some $u$, (cf. [Gu2]).

 \medskip
\noindent (4.3) {\it Any horizontal
geodesic $\alpha$ with $\alpha(0) \in
\partial \Omega_{u}$ must be entirely contained in $
\partial \Omega_{u}$. Consequently, $\Cal H_p \subset T^-_p(
\partial \Omega_{u})$,
where $T^-_p( \partial \Omega_{u})$ is
the tangent cone of $
\partial \Omega_{u}$ at $p$. }
\medskip

Recall that by (2.2) we have
$$\Cal N^+_{p}(\Omega_u, M^n) \subset \big[\Cal
N^+_{p}(\Omega_u, \text{int}(\Omega_{u+ \epsilon}) ) \oplus \Cal
N^+_{p}(\text{int}(\Omega_{u+ \epsilon}), M^n)\big].$$

For $\vec v$ in either  $\Cal N^+_{p}(\Omega_u,
\text{int}(\Omega_{u+ \epsilon}))$ or $ \Cal
N^+_{p}(\text{int}(\Omega_{u+ \epsilon}), M^n) $, we will  show that such a $\vec v$ satisfies (4.2).

    It follows from Theorem 1.10 of [ChG] (or Corollary 1.4 of
    [Y1]) that any {\it  minimal normal vector} $\vec v$ of $\Cal N^+_p(\Omega_u,
     \text{int}(\Omega_{u + \epsilon}))$
    stays {\it minimal} under parallel transport along any
    geodesic in $\partial \Omega_u$.
    Since the convex hull of minimal normal vectors is equal to the outward
normal cone (cf. [Y1, Proposition 1.7]),
the bundle $\Cal N^+(\Omega_u, \text{int}(\Omega_{u + \epsilon}))$ is
    invariant under parallel transport along any  geodesic in $\partial \Omega_u$.
This together with (4.3) implies that if $\vec v \in \Cal N^+_{p}\big(\Omega_u,
\text{int}(\Omega_{u+ \epsilon})\big)$ then $\vec v$ satisfies (4.2).

For $ \vec v \in \Cal N^+_{p}( \text{int}(\Omega_{u+\epsilon}),
M^n)$, the assertion (4.2) follows from Corollary 3.2 of [Gu1]. In
fact, since $\text{int}(\Omega_{u+\epsilon}) $ is totally convex
and totally geodesic,  both $T(\text{int}(\Omega_{u+\epsilon}) )$
and $N^+(\text{int}(\Omega_{u+ \epsilon}), M^n)$ are invariant
under parallel transport along any geodesic in
$\text{int}(\Omega_{u+ \epsilon})$. This together with (4.3)
implies that
 (4.2) holds for any $ \vec v \in \Cal N^+_{p}(
\text{int}(\Omega_{u+\epsilon}), M^n)$.

 Therefore, (4.2) holds for any $ \vec v \in \Cal
N^+_p(\Omega_{u}, M^n)$. This completes the proof of Theorem 4.
\qed\enddemo

    In order to see that  Theorem 4 implies Theorem 1,
 we need to establish a bootstrap argument for the smoothness of ruled
    surfaces. A $C^1$-smooth one-parameter family of a straight
    lines in $\Bbb R^3$ gives rise to a ruled surface. Suppose
    that
    $\{\beta(s), \vec v(s) \}$ are $C^1$-smooth vector valued functions with
    $[\beta'(s)+ t \vec v'(s)] \wedge \vec v(s) \neq 0$ for all $(s, t) \in (a, b)
    \times (c, d)$.
    Then we have a corresponding $C^1$-smooth  immersed ruled surface.
    $$
     \aligned
         F: \quad & (a, b) \times (c, d) \to \Bbb R^3 \\
             &  (s, t) \to \beta(s) + t \vec v(s)
     \endaligned
    $$

    Our bootstrap argument is motivated by the following
    observation.

\medskip
\proclaim{Lemma 5 } (The smooth extension for ruled surfaces in
$\Bbb R^3$) Let $F\big( (a, b) \times (c, d)
\big) = \Sigma^2 $ be an embedded ruled surface in $\Bbb R^3$ and let $F: (a, b) \times (c, d) \to \Bbb R^3$ be a
$C^{1,1}$-smooth embedding map  be as above. Suppose that a subset
$\hat \Sigma^2_\epsilon = F\big( (a, b) \times (\epsilon_1,
\epsilon_2) \big)$ is a $C^\infty$-smooth embedded surface of
$\Bbb R^3$, where $(\epsilon_1, \epsilon_2) \subset (c, d)$. Then
the whole ruled surface $\Sigma^2 $ is a $C^\infty$-smooth surface of $\Bbb R^3$.
\endproclaim
\demo{Proof} By our assumption, $F$ is an embedding map, and hence
the surface $\hat \Sigma^2_\epsilon = F\big( (a, b) \times
(\epsilon_1, \epsilon_2) \big)$ is foliated by straight lines.
Because the surface $\hat \Sigma^2_\epsilon$ and each orbit (each
straight line) are $C^\infty$, the quotient space $Q = [\hat
\Sigma^2 / \sim]$ is a $C^\infty$-smooth 1-dimensional space as
well, where $\sim$ is the equivalent relation induced by the
orbits (the ruling straight lines). Thus, we have a fibration $
(\epsilon_1, \epsilon_2) \longrightarrow \hat \Sigma^2_\epsilon
\longrightarrow Q$. We may assume that $Q = (0, 1)$. Let $\pi:
\hat \Sigma^2_\epsilon \to Q$ be the quotient map. Because the
fibration is topologically trivial, we can find two {\it disjoint}
$C^\infty$-smooth cross-sections
$$
\aligned
 h_i: \quad & Q \to \hat \Sigma^2_\epsilon \\
       &  u \to h_i(u)
       \endaligned
$$
for $i = 0, 1$, where $\pi(h_i(u)) = 0$. (Since the fibre is
1-dimensional line, we may assume that the graph of the
cross-section $h_1$ lies above that of $h_0$ ). Because $h_0(Q)$
and $h_1(Q)$ are disjoint, we obtain a new $C^\infty$-smooth
parametrization of the ruled surface
$$  \aligned
G: \quad & Q \times \Bbb R \to \Bbb R^3 \\
    & (u, \lambda) \to  h_0(u) +  \lambda \frac{[h_1(u) - h_0(u)
    ]}{ \|h_1(u) - h_0(u)   \|   }
\endaligned
$$
Clearly, $G$ is a $C^\infty$-smooth map
 with
$\Sigma^2 \subset G( Q \times \Bbb R )$. Because $F$ is an
embedding map, on the subset $G^{-1}(\Sigma^2)$, one can check
that $G$ remains to be injective and with non-vanishing Jacobi
$G_u \wedge G_\lambda \neq 0$. Hence, $G|_{G^{-1}(\Sigma^2)}$ is
an embedding as well. Thus, $\Sigma^2$ is a $C^\infty$-smooth
embedded surface.
 \qed\enddemo

  The proof of Lemma 5 can be applied to the proof of Theorem 1 as
  follows. Let $\Omega_u$ be a totally convex subset  as above. By
  Federer's Theorem, the hypersurface $\partial[
  U_\delta(\Omega_u)]$ is  $C^{1, 1}$-smooth if $\delta$ is less
  than the cut-radius of $\Omega_u$. Assume that $T > u $ and $d =
  \delta_T - \delta > 0$. Let $\vec v(x)$ be the outward unit normal vector
  of $\partial [U_\delta(\Omega)]$ at $x$.  There is an embedding:
  $$
  \aligned  F: \quad & \partial[U_\delta(\Omega)] \times (c, d) \to M^n
  \\ & (x, t) \to \text{Exp}_x[t \vec v(x)]
\endaligned \tag6.1
  $$
where $c = - \delta$.

\medskip
\proclaim{Proposition 6 } (The smooth extension for the ruled
sub-manifold) For each $y \in \partial[U_\delta(\Omega)]$, we let
$B^{k-1}_r(y) \subset
\partial[U_\delta(\Omega)]$ be a small $k$-dimensional ball around
$y$ which is $C^1$-diffeomorphic to $B^{k-1}(0) \subset \Bbb R^k$
and let $\Omega_u, \delta, \delta_T, c, d$ and $F$ be as above.
Suppose that $F$ is an $C^{1,1}$-smooth embedding and that  $\hat
\Sigma^k_\epsilon = F\big( B^{k-1}_r(y) \times (\epsilon_1,
\epsilon_2) \big)$ is a $C^\infty$-smooth embedded $k$-submanifold
of $M^n$, where $(\epsilon_1, \epsilon_2) \subset (c, d)$. Then
the whole ruled submanifold $\Sigma^k = F\big( B^{k-1}_r(y) \times
(c, d) \big)$ is a $C^\infty$-smooth submanifold of $M^n$.
\endproclaim
\demo{Proof} The proof of Proposition 6 is the same as above with
minor modifications. By our assumption, $\hat \Sigma^k_\epsilon$
is foliated by $C^\infty$-smooth open geodesic segments. The
quotient space $Q = [\hat \Sigma^k_\epsilon/\sim]$ is a
$C^\infty$-smooth $(k-1)$-dimensional open manifold. Because the
fibration $ (\epsilon_1, \epsilon_2) \longrightarrow \hat
\Sigma^k_\epsilon \longrightarrow Q$ is trivial, we can choose two
{\it disjoint} cross sections $h_0: Q \to \hat \Sigma^k_\epsilon$
for $i =0, 1$. If $\pi: \hat \Sigma^k_\epsilon \to Q$ is the
quotient map, then $\pi \circ h_i (u) = u $ for all $u \in Q$.
Since the two cross-sections are disjoint, we may assume that
$r(h_1(u)) > r(h_0(u))$ for all $u \in Q$, where $r(y) = d\big(y,
\partial[U_\delta(\Omega)] \big)$. For each $u \in Q$, we consider
the unit vector
$$\vec \eta(u) =
\frac{\text{Exp}_{h_0(u)}^{-1}[h_1(u)]}{\|
\text{Exp}_{h_0(u)}^{-1}[h_1(u)]  \|}
$$
at the point $h_0(u)$. Similarly, we consider a new
$C^\infty$-smooth parametrization
$$  \aligned
G: \quad & Q \times \Bbb R \to M^n \\
    & (u, \lambda) \to  \text{Exp}_{h_0(u)}[ \lambda \vec \eta (u)
    ].
\endaligned
$$
Clearly, we have $ \Sigma^k = F\big( B^{k-1}_r(y) \times (c, d)
\big) \subset G\big(Q \times \Bbb R\big)$. This completes the
proof. \qed\enddemo

   With Lemma 2, Theorem 4 and Proposition 6, we are ready to prove Theorem 1.

\bigskip
\demo{Proof of Theorem 1} Let $\{ \Omega_u \}$ be a
Cheeger-Gromoll convex exhaustion described as above.  It is
sufficient to verify that the subset
   $[U_{\delta_0(T)}(\Omega_T) \cap F_q]$ has a $k$-dimensional
   $C^\infty$-smooth interior for any given $T > a_m$ and $q \in \Cal S$, where $\delta_0(T)$ is given by
   Lemma 2.

   Fix $T > a_m$ with $\dim[\Omega_T] = n$. Let $\delta_0(T)$ be
   given by Lemma 2 and $C_T$ be given by (2.3). Choose a
   partition $0=u_0 < u_1 < ... < u_N = T$ of $[0, T]$ such that
   $u_{j} - u_{j-1} < \frac{2C_T}{\delta_0(T)}$ for $j=1, ...,N,$
   where $N= N_T$ is a number depending on $T$.

 We will prove the following assertion by induction
on $j=0, 1, ..., N$.
   \medskip
   \noindent
    {\bf Assertion j.} {\it The sub-level set
    $[U_{\delta_0(T)}(\Omega_{u_j}) \cap F_q]$ has the $k$-dimensional
   $C^\infty$-smooth interior, where   $q \in \Cal
   S$ and $k = \dim[F_q]$.}

   \medskip

It follows from Perelman's theorem or Theorem
   4 that $\text{Exp}_q[\Cal N_q^+(\Cal S, M^n)      ] \subset
   F_q$.
   Since the soul $\Cal S$ has the cut radius $\ge \delta_0(T)$
   and $\Cal S$ is $C^\infty$-smooth, Assertion 0 holds.

Let $ \epsilon_1 = \frac{\delta_0(T)}{16}  $  and $\epsilon_2 =
\frac{\delta_0(T)}{8}$. We consider
$$ 
         A(\Omega_{u_j}, r_1, r_2) = \{ z \in F_q | \quad 0 < r_1 <  d(z,
         \Omega_{u_j}) < r_2 \}
$$
It is clear that  $A(\Omega_{u_1}, \epsilon_1, \epsilon_2) \subset
U_{\delta_0(T)}(\Cal S)$. It follows from Assertion 0 that  the
subset $\hat \Sigma^k_\epsilon =  A(\Omega_{u_1}, \epsilon_1,
\epsilon_2) \subset F_q \cap U_{\delta_0(T)}(\Cal S)$ is
$C^\infty$-smooth $k$-dimensional open sub-manifold. By Theorem 4,
we let $\Sigma^k_1 = A(\Omega_{u_1},\frac{\delta_0(T)}{16},
\delta_0(T) )$ be the ruled $k$-dimensional submanifold. It
follows from Proposition 6 (the smooth extension theorem for the
ruled submanifold) that $\Sigma^k_1$ is a $C^\infty$-smooth
$k$-dimensional submanifold of $M^n$. Observe that the subset
$[U_{\delta_0(T)}(\Omega_{u_j}) \cap F_q] $ is contained in the
union $ \{ [U_{\delta_0(T)}(\Cal S) \cap F_q] \cup \Sigma^k_1 \}$.
Since $\Sigma^k_1$ is a $C^\infty$-smooth, Assertion 1 follows
from Assertion 0.

            Similarly, using Theorem 4 and Proposition 6 we can  verify
            that Assertion (j-1) is true then
            Assertion j holds as well for $j\ge 2$. In fact, by induction    we
see that  $ A(\Omega_{u_j}, \epsilon_1, \epsilon_2) \subset
[U_{\delta_0(T)}(\Omega_{u_{j-1}}) \cap F_q]$ is
$C^\infty$-smooth. It follows from Theorem 4 and Proposition 6
that the ruled submanifold $\Sigma^k_j =
A(\Omega_{u_j},\frac{\delta_0(T)}{16}, \delta_0(T) )$ must be of
$C^\infty$-smooth as well. Since $ [U_{\delta_0(T)}(\Omega_{u_j})
\cap F_q] \subset [U_{\delta_0(T)}(\Omega_{u_{j-1}}) \cap F_q]
\cup \Sigma^k_j$, Assertion j follows. Theorem 1 follows from
Assertion $N_T$ for any arbitrarily large T. \qed\enddemo

     The first author is very grateful to Professor  B. Wilking for pointing out a mistake in
    an earlier version of  the manuscript [CaS].

\enddocument